\documentclass[preprint,12pt]{elsarticle}
\usepackage[latin9]{inputenc}
\usepackage{color}
\usepackage{amsthm}
\usepackage{amsmath}
\usepackage{amssymb}
\usepackage{graphicx}
\usepackage{esint}

\journal{Computers and Mathematics with Applications}

\makeatletter
\theoremstyle{plain}
\newtheorem{thm}{\protect\theoremname}
  \theoremstyle{plain}
  \newtheorem{lem}[thm]{\protect\lemmaname}

\@ifundefined{date}{}{\date{}}

\usepackage{amsthm}

\newcommand{\e}{\mathrm{e}}

\makeatother

  \providecommand{\lemmaname}{Lemma}
\providecommand{\theoremname}{Theorem}

\begin{document}

\begin{frontmatter}
	
	\title{An almost symmetric Strang splitting scheme for nonlinear evolution equations\tnoteref{label1}}
	\tnotetext[label1]{This work is supported by the Fonds zur F\"orderung der Wissenschaften (FWF) -- project id: P25346.}

	\author[uibk]{Lukas Einkemmer\corref{cor1}}
	\ead{lukas.einkemmer@uibk.ac.at}

	\author[uibk]{Alexander Ostermann}
	\ead{alexander.ostermann@uibk.ac.at}

	\address[uibk]{Department of Mathematics, University of Innsbruck, Austria}
	\cortext[cor1]{Corresponding author}



\begin{abstract}
In this paper we consider splitting methods for the time integration
of parabolic and certain classes of hyperbolic partial differential
equations, where one partial flow can not be computed exactly. Instead,
we use a numerical approximation based on the linearization of the vector field. This is of interest
in applications as it allows us to apply splitting methods to a 
wider class of problems from the sciences.

However, in the situation described the classic Strang splitting scheme,
while still a method of second order, is not longer symmetric. This,
in turn, implies that the construction of higher order methods by
composition is limited to order three only. To remedy this situation,
based on previous work in the context of ordinary differential equations,
we construct a class of Strang splitting schemes that are symmetric
up to a desired order.

We show rigorously that, under suitable assumptions on the nonlinearity,
these methods are of second order and can then be used to construct higher
order methods by composition. In addition, we illustrate the theoretical
results by conducting numerical experiments for the Brusselator system
and the KdV equation.
\end{abstract}

\end{frontmatter}

\section{Introduction\label{sec:Introduction}}

In this paper we consider the time discretization of parabolic and
certain classes of hyperbolic partial differential equations. More
specifically, we assume that the considered problem 
can be written as the following abstract Cauchy problem 
\begin{equation}
u^{\prime}  = Au+B(u),\qquad
u(0)  =u_{0},
\end{equation}
where $A$ is a linear (but possibly unbounded) operator. In this
context splitting methods can be applied if the partial flows generated
by $A$ and $B$ have an analytical representation, or if an efficient algorithm
for finding their exact solution is known. For a review of splitting methods
we refer the reader to \cite{mclachlan1995}. However, the assumption
that the partial flow generated by the nonlinear operator $B$ can
be solved exactly is usually a very strong requirement. Most partial
differential equations which are drawn from the sciences do neither
admit an analytical solution nor can they be solved exactly in an
efficient manner by some algorithm from the literature.

To remedy this deficiency of classic splitting methods, we propose and analyze splitting schemes which approximate
the partial flow generated by $B$, i.e.
\begin{equation}
u^{\prime}=B(u),\label{eq:partial-flow-B}
\end{equation}

by that of an inhomogeneous linear differential equation. That is, instead
of equation \eqref{eq:partial-flow-B}, we consider the linearized
problem given by
\begin{equation}
u^{\prime}=b(u_{\star})u+d,\label{eq:linearized-flow}
\end{equation}
where for consistency reasons we have to assume that 
\[
B(u)=b(u)u+d.
\]

That is, once a value $u_{\star}$ is substituted, the flow corresponding
to \eqref{eq:linearized-flow} can be computed efficiently. That such
a linearized problem can be solved exactly is a much less stringent
condition and thus splitting methods can potentially be applied to
a larger class of problems.

For the remainder of the paper let us denote the flow corresponding
to equation \eqref{eq:linearized-flow} by $\varphi_{t}^{b(u_{\star})}$,
which can be written explicitly, by employing the exponential
and $\phi_{1}$ function, as follows
\begin{equation}
\varphi_{t}^{b(u_{\star})}(u(0))=\e^{tb(u_{\star})}u(0)+t\phi_{1}\left(tb(u_{\star})\right)d,
\end{equation}
where the $\phi_{1}$ function is given by
\begin{equation}
\phi_{1}(z)=\frac{\e^{z}-1}{z}.
\end{equation}
In this context we can formulate the (classic) Strang splitting scheme
as follows
\begin{equation} \label{eq:classic-strang}
\begin{aligned}
u_{1/2} & =\varphi_{\frac{\tau}{2}}^{b(u_{0})}(\e^{\frac{\tau}{2}A}u_{0}),\\
u_{1}=M_{\tau}(u_{0}) & =\e^{\frac{\tau}{2}A}\varphi_{\tau}^{b(u_{1/2})}(\e^{\frac{\tau}{2}A}u_{0}). 
\end{aligned}
\end{equation}
Such a splitting scheme has been employed, for example, to numerically
solve the Vlasov--Poisson equations (see e.g. \cite{Cheng:1976}
and \cite{time_analysis}) or the Davey--Stewartson equation (see
e.g. \cite{klein2011}). Note, however, that in both of these cases
the underlying structure of the equations guarantees that the Strang
splitting scheme, as stated above, is still symmetric. Therefore,
the usual construction of composition methods (as described e.g. in
\cite{hairer2006}) can be carried out without modification. 

However, for many problems arising in the sciences, such as nonlinear
reaction-diffusion equations or the KdV equation (both of which are
discussed in section \ref{sec:Numerical-simulations}), such a simplification
can not be assumed. Then the Strang splitting scheme above is, unfortunately,
no longer symmetric. The lacking symmetry of the method does not severely
affect performance; however, if composition is used as a means to
construct higher order methods, symmetry is usually a necessary condition.

To remedy this situation a family of Strang
splitting schemes were introduced in the context of ordinary differential
equation (see \cite{einkemmer2013}). The construction of these schemes is based on the fact that the Lie splitting algorithm 
\begin{equation}
L_{\frac{\tau}{2}}=\varphi_{\frac{\tau}{2}}^{b(u_{0})}\circ\e^{\frac{\tau}{2}A}\label{eq:Lie-half-step}
\end{equation}
composed with its adjoint method, which we denote by $L_{\frac{\tau}{2}}^{*}$, is a method of second order and symmetric. It is given by
\[
S_{\tau}=L_{\frac{\tau}{2}}^{*}\circ L_{\frac{\tau}{2}}.
\]
However, for our purpose it it is more convenient to represent $u_1=S_{\tau}u_0$ as the solution of the following implicit equation
\[
u_{1}=\e^{\frac{\tau}{2}A}\circ\varphi_{\frac{\tau}{2}}^{b(u_{1})}(u_{1/2}).
\]
Note that to solve this equation is computationally not attractive in practice. Therefore, it has
been suggested to employ $i\in\mathbb{N}$ fixed-point iterations
to compute an approximation to $u_{1}$. Let us denote the resulting
scheme by $S_{\tau}^{(i)}$. The corresponding starting value (for
the fixed-point iteration) is given by
\begin{align} \label{eq:Stau1}
S_{\tau}^{(1)} & =\e^{\frac{\tau}{2}A}\circ\varphi_{\frac{\tau}{2}}^{b(u_{1/2})}\circ\varphi_{\frac{\tau}{2}}^{b(u_{0})}\circ \e^{\frac{\tau}{2}A}.
\end{align}
The one-step methods $S_{\tau}^{(i)}$ are not symmetric but they
are symmetric up to order $i$; that is, they satisfy the following
relation $\big(S_{\tau}^{(i)}\big)^{*}=S_{\tau}^{(i)}+\mathcal{O}\left(\tau^{i}\right)$.
Therefore, it is possible to construct composition methods of arbitrary
(even) order, where $i$, and thus the computational cost of $S_{\tau}^{(i)}$,
increases linearly with the desired order $p$. For a more detailed
description we refer the reader to \cite{einkemmer2013}.

In this paper, we will provide a rigorous convergence analysis which
shows that for parabolic and certain classes of hyperbolic partial
differential equations the method $S_{\tau}^{(i)}$ is of second order
for $i\geq2$. This analysis conducted in the context of partial differential
equations is significantly more involved than the analysis conducted
in \cite{einkemmer2013}. The assumptions we make in order to prove
the main result will be that the classic Strang splitting scheme given in equation \eqref{eq:classic-strang} is convergent of order two
(i.e. that the classical Strang splitting scheme can be successfully applied to the problem under consideration) and that terms of the form
\[
b^{i}(\tilde{u})b^{(\ell)}(\tilde{u})b^{j}(\tilde{u})u(t),\qquad \ell+i+j\leq3
\]
are bounded in a suitable norm; here $\tilde{u}$ is, in a sense to
be made precise, close to the exact solution $u(t)$. Note that in
the linear case this assumption is equivalent to the fact that $b^{3}u(t)$
can be bounded; i.e.~we have to bound three applications of the operator
$b$ to the exact solution. It is useful to keep that in mind in order
to compare the assumptions made here to, for example, \cite{jahnke2000},
where the analysis of the exact Strang splitting scheme in the linear
case is provided.

To that end we will discuss the fixed-point iteration in the stiff
case (section \ref{sec:Fixed-point-iteration}), and show second order convergence of $S_{\tau}^{(i)}$
for the hyperbolic (section \ref{sec:Convergence-for-hyperbolic})
and parabolic (section \ref{sec:Convergence-for-parabolic}) case.
In section \ref{sec:high-order-composition} we discuss the construction
of higher order composition methods. Finally, we will show in section
\ref{sec:Numerical-simulations} that in the case of the KdV equation
as well as the Brusselator system the predicted behavior agrees with
the numerical experiments conducted.

\section{Fixed-point iteration in the stiff case\label{sec:Fixed-point-iteration}}

The purpose of this section is to discuss the fixed-point iteration
necessary to compute $S_{\tau}^{(i)}$, as outlined in section \ref{sec:Introduction},
in more detail. In the non-stiff case a sufficiently small $\tau$
can always be found such that the iteration converges. However, in
the present discussion it is necessary to employ a condition on the,
in general, unbounded nonlinearity $b$. We will assume that a Lipschitz
condition is satisfied in a Banach space $\left(V,\Vert\cdot\Vert\right)$
that includes, in a sense to me made precise in section \ref{sec:Convergence-for-hyperbolic},
functions close to the exact solution.
\begin{thm}\label{thm:fixedpoint}
Suppose that $b(\cdot)$ is Lipschitz continuous in $V$, i.e. there
exists a constant $L$ such that 
\[
\Vert b(v_{1})w-b(v_{2})w\Vert\leq L\Vert v_{1}-v_{2}\Vert
\]
for all $v_{1},v_{2},w\in V$. 
Then there exists a step size $\tau_0>0$, which depends only on the Lipschitz contant $L$ and on bounds of the partial flows, such that for all \mbox{$\tau \in [0,\tau_0)$} the fixed-point iteration convergences. In addition, it holds that 
\[
\Vert S_{\tau}^{(i)}-S_{\tau}\Vert\leq C\left(CL\right)^{i-1}\tau^{i+1},\quad i\geq1.
\]
\end{thm}
\begin{proof}
We have to consider the fixed-point iteration given by
\begin{equation} \label{eq:fixedpointiteration}
u = F(u)=\e^{\frac{\tau}{2}A}\varphi_{\frac{\tau}{2}}^{b(u)}(u_{1/2}).
\end{equation}
By employing the variation-of-constants formula, we can write 
\[
	\varphi_{\frac{\tau}{2}}^{b(v_{1})}(u_{1/2})=\varphi_{\frac{\tau}{2}}^{b(v_{2})}(u_{1/2})+\int_{0}^{\frac{\tau}{2}}\e^{(\frac{\tau}{2}-s)b(v_{2})}\left(b(v_{1})-b(v_{2})\right)\varphi_{s}^{b(v_{1})}(u_{1/2})\,\mathrm{d}s
\]
and by using the Lipschitz continuity of $b$, we get 
\[
\Vert F(v_{1})-F(v_{2})\Vert\leq CL\tau\Vert v_{1}-v_{2}\Vert.
\]
Thus, our fixed-point iteration convergences, if $\tau_0$, as given in the statement of the theorem, satisfies $CL\tau_0\leq1$.
\end{proof}
It should be duly noted that although the step size is limited by
the above theorem this is not a CFL condition. In our case the step
size does only depend on properties of $b(u)$, the exact solution as
well as the scheme under consideration. However, the bound holds independently
of the specific space discretization under consideration. In fact,
the above result is what one would expect for a stable explicit scheme;
namely that the step size is limited by the Lipschitz constant of
the right-hand side.

In the next section we will employ heavily the following result; it
states that a time step conducted by applying $S_{\tau}$ to some
$u_{0}$ results in a value that is different from $u_{0}$ by at
most an order of $\tau$. For the convenience of the reader we state
and proof this well-known result in the next theorem.
\begin{lem} 
\label{lem:first-order-bound} If the Lie splitting algorithm \eqref{eq:Lie-half-step}
is consistent of order $1$ and $\tau$ is chosen such that the fixed-point iteration \eqref{eq:fixedpointiteration} converges,
then 
\[
\Vert S_{\tau}(u_{0})-u_{0}\Vert\leq C\tau,
\]
where the constant $C$ is independent of $\tau$.\end{lem}
\begin{proof}
From equation \eqref{eq:Stau1} we know that $S_{\tau}^{(1)}$ is the composition of two Lie steps. Therefore, it holds that
\[
	\Vert S_{\tau}^{(1)}(u_0)-u_0 \Vert \leq C_2 \tau.
\]
By employing the uniqueness and existence result for the fixed-point
iteration, we get 
\[
\Vert S_{\tau}^{(i)}(u_{0})-u_{0}\Vert\leq\Vert S_{\tau}^{(1)}(u_{0})-u_{0}\Vert\sum_{k=0}^{i-1}\left(C_1L\tau\right)^{k}\leq\frac{C_2}{1-C_1L\tau}\tau,
\]
for any $i\geq2$. This bound gives the desired result, since under the assumptions of convergence of the fixed-point iteration we know that $C_1L\tau<1$ (see Theorem \ref{thm:fixedpoint}).
\end{proof}

\section{Convergence for hyperbolic differential equations\label{sec:Convergence-for-hyperbolic}}

First, we discuss the condition
on the initial value that is necessary to achieve a given regularity in
time. This will give us an idea of the required regularity assumptions on the initial value, which is needed for the schemes described here to be convergent of second order. To accomplish this, in section \ref{sub:Regularity-in-time}, a suitable toy model is analyzed.
In section \ref{sub:Convergence-of-the-Strang} we will derive,
for the hyperbolic case, the conditions under which the splitting
scheme $S_{\tau}^{(i)}$, for $i\geq2$, is convergent of order two.
We will formulate a condition that depends on the regularity of the
exact solution. More formally, we have to bound a number of applications
of $b(u_{\star})$ and its derivative to the exact solution and assume
that the (classic) Strang splitting scheme $M_{\tau}$ is convergent
of order two. For hyperbolic systems which do not change the regularity
of the initial value this is equivalent to assuming some regularity
for the initial value.

\subsection{Regularity in time\label{sub:Regularity-in-time}}

An interesting question is how much regularity on the initial
value is required such that a scheme can be convergent of second order in time. To elaborate
on this question let us consider the following (nonlinear) advection
equation 
\begin{equation}
\partial_{t}u(t,x)+u(t,0)\partial_{x}u(t,x)=0.\label{eq:toy_model}
\end{equation}
The solution of this toy model can be written as 
\begin{equation}
u(t,x)=u_{0}\left(x-v(t)\right),\label{eq:toy_exact}
\end{equation}
where the advection speed depends on the solution as follows 
\[
v(t)=\int_{0}^{t}u(\sigma,0)\,\mathrm{d}\sigma.
\]
Now, let us investigate the time regularity of the solution $u(t,x)$.
To that end let us assume that for the initial value $u_{0}$ it is
true that $u_{0}\in\mathcal{C}^{1,\alpha}$ for some $\alpha\in(0,1)$.
From \eqref{eq:toy_model} we deduce that $u(t,x)$ is at least once
differentiable in time. The resulting function $\partial_{t}u(t,\cdot)$
lies in $\mathcal{C}^{0,\alpha}$; this can be easily deduced from
equation \eqref{eq:toy_exact}. Our goal is to show that 
\[
\frac{\vert\partial_{t}u(t,x)-\partial_{t}u(s,x)\vert}{\vert t-s\vert^{\beta}}
\]
is not bounded for any $\beta>\alpha$ as $s$ tends to $t$. Rewriting
the above expression we get 
\begin{align*}
(\partial_{t}u)(t,x)-(\partial_{t}u)(s,x)&=\left(u(t,0)-u(s,0)\right)\partial_{x}u(t,x) \\&\quad+u(s,0)\left((\partial_{x}u)(t,x)-(\partial_{x}u)(s,x)\right).
\end{align*}
The first term is Lipschitz continuous in time and therefore can be
bounded. However, for the second term we can write 
\begin{align*}
\frac{\vert(\partial_{x}u)(t,x)-(\partial_{x}u)(s,x)\vert}{\vert t-s\vert^{\beta}} & =\frac{\vert u_{0}^{\prime}(x-v(t))-u_{0}^{\prime}(x-v(s))\vert}{\vert t-s\vert^{\beta}},
\end{align*}
where by defining $z=x-v(t)$ and $w=x-v(s)$ and using the estimate
\[
\vert z-w\vert=\left\vert \int_{s}^{t}u(\sigma,0)\,\mathrm{d}\sigma\right\vert \geq\min_{\sigma\in[0,t]}\vert u(\sigma,0)\vert\vert t-s\vert
\]
we get 
\[
\frac{\vert u_{0}^{\prime}(x-v(t))-u_{0}^{\prime}(x-v(s))\vert}{\vert t-s\vert^{\beta}}\geq\frac{\vert u_{0}^{\prime}(z)-u_{0}^{\prime}(w)\vert\cdot\min_{\sigma\in[0,t]}\vert u(\sigma,0)\vert^{\beta}}{\vert z-w\vert^{\beta}}.
\]

If we now choose a $u_{0}\in\mathcal{C}^{1,\alpha}$ that nowhere
can be improved to a better Hölder exponent $\beta>\alpha$ and is
bounded away from zero (which then also holds true for $u(t,\cdot)$)
the term on the right hand side diverges as $z$ tends to $w$ and
thus we have shown that the regularity of $u(t,x)$ in time is at
most $\mathcal{C}^{1,\alpha}$. That such an initial value can be
found is shown in \cite{katzourakis2010}.

As a numerical approximation to the solution of \eqref{eq:toy_model}, we consider the exact solution of
\begin{equation} \label{eq:toy_model_simplified}
\partial_{t}u(t,x)=u\left(\tfrac{\tau}{2},0\right)\partial_{x}u(t,x)
\end{equation}
on the interval $[0,\tau]$. This is not a practical numerical method
as we have, in general, no algorithm that allows us to compute $v(\tau/2)$
exactly. However, it can be seen as the limit of a scheme, where we
approximate $v$ by approximating $u(t,0)$ in some suitable manner.
In this case the numerical solution, i.e. the exact solution of equation \eqref{eq:toy_model_simplified}, is given by the midpoint rule, i.e.
\[
u_{1}=u_{0}\left(x-\tau u(\tfrac{\tau}{2},0)\right),
\]
which has to be compared to the exact solution, given in \eqref{eq:toy_exact}.
By using the method of characteristics we can estimate the error as
follows 
\[
\Vert u(\tau,x)-u_{1}\Vert\leq\tau L\left\Vert u(\tfrac{\tau}{2},0)-\frac{1}{\tau}\int_{0}^{\tau}u(t,0)\mathrm{d}t\right\Vert ,
\]
where $L$ is the Lipschitz constant of $u_{0}$. Therefore, the numerical
scheme can be seen as the midpoint rule approximating an integral.
It is well-known that the midpoint rule yields an approximation with
a local error of size $C\tau^{\min(3,k+\alpha)}$, if and only if the integrand lies in $\mathcal{C}^{k,\alpha}$.

\subsection{Convergence of the iterated Strang splitting scheme\label{sub:Convergence-of-the-Strang}}

The discussion in the previous section shows that we can expect the
Strang splitting scheme to be second order accurate in time only if
we require that $u_{0}\in\mathcal{C}^{2,1}$. In the framework discussed
in section \ref{sec:Introduction} this means that we have to bound
at least three applications of the operator $b(u_{\star})$ to the
exact solution. In fact our assumption, which we will state in this
section, will turn out to be slightly more complicated due to the
fact that we have to bound certain derivatives of $b(\cdot)$ as well.

Before stating the assumption that is used to prove the consistency
result, let us remind the reader that for a function $f\colon\, V\times W\to Z$
which is Fr\'echet differentiable with respect to the first component
it holds that 
\[
\Vert(\partial_{1}f)(u_{0},w)d\Vert\leq\Vert(\partial_{1}f)(u_{0},w)\Vert\Vert d\Vert,
\]
where it is understood that $d\in V$. Thus, we can separate the condition
on the function $f$ from the direction of the derivative. We will
see below that this fact simplifies the condition needed to establish
the correct order of the scheme significantly. Note that, clearly,
this also holds true if the function $f$ is linear in the second
component and we employ the notation outlined in section \ref{sec:Introduction}. 

Now, let us assume the following bound
\begin{equation}
	\sup_{t\in[0,T]} \left\Vert b^{i}(\tilde{u})b^{(\ell)}(\tilde{u})b^{j}(\tilde{u})u(t)\right\Vert \leq C
	,\qquad \ell+i+j\leq3,\label{eq:assumption_1}
\end{equation}
where $\ell\in\left\{ 0,1,2\right\} $, $i,j\in\mathbb{N}$.
The function $\tilde{u}$ are generic, i.e.~different occurrences of such functions
do not need to be indistinguishable, linear combinations of the exact
solution $u(s)$ as well as $L_{\frac{\tau}{2}}\left(u(s)\right)$
and $S_{\tau}\left(u(s)\right)$ for some $s\in\left[0,T\right]$. 

The following theorem shows consistency under the additional assumption
that the (classic) Strang splitting scheme $M_{\tau}$ is of the expected
order.
\begin{thm}
\label{thm:(Consistency)}(Consistency) Let us consider the following
abstract initial value problem 
\begin{equation*}
u^{\prime}  =Au+b(u)u+d,\quad
u(0)  =u_{0}
\end{equation*}
 and suppose that assumption \eqref{eq:assumption_1} is satisfied.
If, in addition, the Strang splitting scheme $M_{\tau}$ is consistent
of order two, we can conclude that $S_{\tau}^{(i)}$, for $i\geq2$,
is consistent of order two, i.e. 
\[
\Vert S_{\tau}^{(i)}u_{0}-u(\tau)\Vert\leq C \tau^{3}.
\]
\end{thm}
\begin{proof}
First, let us consider a single time step of the implicit scheme $S_{\tau}$, i.e. $u_{1}=S_{\tau}u_{0}$,
can be written as (note that $u_{0}$ is a not to be confused with
the initial value; it is the exact solution at a point in time to
which a splitting operator is applied) 
\[
u_{1}=\varphi_{\frac{\tau}{2}}^{A}\circ\varphi_{\frac{\tau}{2}}^{b(u_{1})}\circ\varphi_{\frac{\tau}{2}}^{b(u_{0})}\circ\varphi_{\frac{\tau}{2}}^{A}(u_{0}).
\]
Now let us proceed by comparing our method, i.e. $S_{\tau}$, to the
Strang splitting scheme $M_{\tau}$. The difference can be estimated
as follows 
\begin{align*}
&\left\Vert \varphi_{\frac{\tau}{2}}^{A}\circ\varphi_{\tau}^{b(u_{1/2})}\circ\varphi_{\frac{\tau}{2}}^{A}(u_{0})-\varphi_{\frac{\tau}{2}}^{A}\circ\varphi_{\frac{\tau}{2}}^{b(u_{1})}\circ\varphi_{\frac{\tau}{2}}^{b(u_{0})}\circ\varphi_{\frac{\tau}{2}}^{A}(u_{0})\right\Vert \\
&\quad\leq\left\Vert \varphi_{\frac{\tau}{2}}^{A}\right\Vert \left\Vert \varphi_{\tau}^{b(u_{1/2})}(v_{0})-\varphi_{\frac{\tau}{2}}^{b(u_{1})}\circ\varphi_{\frac{\tau}{2}}^{b(u_{0})}(v_{0})\right\Vert 
\end{align*}
with $v_{0}=\varphi_{\frac{\tau}{2}}^{A}(u_{0})$. Note that the flow
$\varphi_{\tau}^{b(u_{\star})}$ can be written as 
\[
\varphi_{\tau}^{b(u_{\star})}(u_{0})=\e^{\tau b(u_{\star})}u_{0}+\tau\phi_{1}\left(\tau b(u_{\star})\right)d.
\]
Therefore we have to compare 
\begin{equation}
\varphi_{\frac{\tau}{2}}^{b(u_{1})}\circ\varphi_{\frac{\tau}{2}}^{b(u_{0})}(v_{0})=\e^{\frac{\tau}{2}b(u_{1})}\e^{\frac{\tau}{2}b(u_{0})}v_{0}+\tfrac{\tau}{2}\left[\e^{\frac{\tau}{2}b(u_{1})}\phi_{1}\left(\tfrac{\tau}{2}b(u_{0}))\right)+\phi_{1}\left(\tfrac{\tau}{2}b(u_{1})\right)\right]d\label{eq:S_tau_compare}
\end{equation}
to the following expression 
\begin{equation}
\varphi_{\tau}^{b(u_{1/2})}(v_{0})=\e^{\tau b(u_{1/2})}v_{0}+\tau\phi_{1}\left(\tau b(u_{1/2})\right)d.\label{eq:M_tau_compare}
\end{equation}

 Let us first compare the homogeneous parts. By expanding this part
in equation \eqref{eq:M_tau_compare} 
\[
\e^{\tau b(u_{1/2})}=I+\tau b(u_{1/2})+\tfrac{\tau^{2}}{2}b^{2}(u_{1/2})+\tau^{3}b^{3}(u_{1/2})\phi_{3}\left(\tau b(u_{1/2})\right),
\]
where the $\phi_{k}$ functions are defined recursively by 
\begin{align*}
\phi_{k}(\tau v) & =\frac{1}{k!}+\tau v\phi_{k+1}\left(\tau v\right)\\
\phi_{0}(\tau v) & =\e^{\tau v}.
\end{align*}
To expand equation \eqref{eq:S_tau_compare} let us first define $g(\tau)=\e^{\frac{\tau}{2}b(u_{1})}\e^{\frac{\tau}{2}b(u_{0})}$.
Then 
\begin{align*}
g^{\prime}(\tau) & =\frac{1}{2}\left(b(u_{1})g(\tau)+g(\tau)b(u_{0})\right)\\
g^{\prime\prime}(0) & =\frac{1}{4}\left(b^{2}(u_{1})+2b(u_{1})b(u_{0})+b^{2}(u_{0})\right)
\end{align*}
which gives the desired expansion 
\[
\e^{\frac{\tau}{2}b(u_{1})}\e^{\frac{\tau}{2}b(u_{0})}=I+\frac{\tau}{2}\left[b(u_{0})+b(u_{1})\right]+\frac{\tau^{2}}{8}\left[b^{2}(u_{1})+2b(u_{1})b(u_{0})+b^{2}(u_{0})\right]+R_{3},
\]
where 
\[
R_{3}=\tau^{3}\sum_{i=0}^{3}c_{i}\int_{0}^{1}(1-\theta)^{2}b^{i}(u_{1})\e^{\frac{\tau}{2}b(u_{1})}\e^{\frac{\tau}{2}b(u_{0})}b^{3-i}(u_{0})\,\mathrm{d}\theta
\]
for appropriately chosen constants $c_{i}$. Now we have to compare
the corresponding terms of equal order.

\emph{Terms of first order.} By employing the fundamental theorem of
calculus we can show that 
\begin{align*}
& b\left(u_{0}+\theta(u_{1/2}-u_{0})\right)-b(u_{0}) \\
&\quad =\int_{0}^{\theta}\frac{\mathrm{d}}{\mathrm{d}\eta}b\left(u_{0}+\eta(u_{1/2}-u_{0})\right)\,\mathrm{d}\eta\\
 &\quad =\int_{0}^{\theta}b^{\prime}\left(u_{0}+\eta(u_{1/2}-u_{0})\right)(u_{1/2}-u_{0})\,\mathrm{d}\eta\\
 &\quad =\theta b^{\prime}\left(u_{0}\right)(u_{1/2}-u_0)
  +\int_{0}^{\theta}\left(\theta-\eta\right)b^{\prime\prime}\left(u_{0}+\eta v_{1/2}\right)(v_{1/2},v_{1/2})\,\mathrm{d}\eta,
\end{align*}
where $v_{1/2} = u_{1/2}-u_{0}$. For $\theta=1$ this gives the following operator identity
\[
	b\left(u_{1/2}\right)=b(u_{0})+b^{\prime}\left(u_{0}\right)(u_{1/2}-u_{0})+\int_{0}^{1}\left(1-\eta\right)b^{\prime\prime}\left(u_{0}+\eta v_{1/2}\right)(v_{1/2},v_{1/2})\,\mathrm{d}\eta.
\]
The same expansion can be carried out for $b(u_{1})$. Finally, we
have to bound 
\[
b(u_{1/2})-\frac{1}{2}\left(b(u_{0})+b(u_{1})\right)=b^{\prime}(u_{0})\left(u_{1/2}-\tfrac{1}{2}(u_{0}+u_{1})\right)+R_{2}^{b(u_{1/2})}-\tfrac{1}{2}R_{2}^{b(u_{1})},
\]
where for the sake of brevity the remainder terms are denoted by $R_{2}^{b(u_{1/2})}$
and $R_{2}^{b(u_{1})}$, respectively. All the terms can be bounded
by assumption \eqref{eq:assumption_1}.  

\emph{Terms of second order.} We have to estimate 
\begin{align*}
 & b^{2}(u_{1})+2b(u_{1})b(u_{0})+b^{2}(u_{0})-4b^{2}(u_{1/2})\\
 & \quad=\left(b^{2}(u_{1})-b^{2}(u_{1/2})\right)+\left(b^{2}(u_{0})-b^{2}(u_{1/2})\right)+2\left(b(u_{1})b(u_{0})-b^{2}(u_{1/2})\right).
\end{align*}
A bound for the first and second term can easily be found. The third
term can be rewritten as
\[
2\left(b(u_{1})-b(u_{1/2})\right)b(u_{0})+2b(u_{1/2})\left(b(u_{0})-b(u_{1/2})\right)
\]
which together with the estimate 
\[
\left\Vert \left(b(u_{0})-b(u_{1/2})\right)u_{0}\right\Vert \leq\sup_{\eta\in[0,1]}\Vert b^{\prime}\left(u_{0}+\eta(u_{1/2}-u_{0})\right)u_{0}\Vert\Vert u_{0}-u_{1/2}\Vert,
\]
Lemma \ref{lem:first-order-bound} as well as assumption \eqref{eq:assumption_1}
is sufficient to show the desired bound.

\emph{Terms of third order.} We have to bound $R_{3}$ as well as $b^{3}(u_{1/2})\phi_{3}(\tau b(u_{1/2}))$.
The first bound is immediate, for the second we note that since the commutator $\left[b^{3}(u_{1/2}),\phi_{3}(\tau b(u_{1/2}))\right]$ vanishes,
the desired result follows easily. 

Finally, let us compare the two inhomogeneous terms which are given
by 
\begin{align*}
\phi_{1}\left(\tau b(u_{1/2})\right) & =I+\tau b\left(u_{1/2}\right)\phi_{2}\left(\tau b(u_{1/2})\right)\\
 & =I+\tfrac{\tau}{2}b\left(u_{1/2}\right)+\tau^{2}b^{2}\left(u_{1/2}\right)\phi_{3}\left(\tau b(u_{1/2})\right)
\end{align*}
and
\begin{align*}
 & \tfrac{1}{2}\left[\e^{\frac{\tau}{2}b(u_{1})}\phi_{1}\left(\tfrac{\tau}{2}b(u_{0})\right)+\phi_{1}\left(\tfrac{\tau}{2}b(u_{1})\right)\right]\\
 & \qquad =\tfrac{1}{2}\left[I+\e^{\frac{\tau}{2}b(u_{1})}\right]+\tfrac{\tau}{8}\left[\e^{\frac{\tau}{2}b(u_{1})}b(u_{0})+b(u_{1})\right]\\
 & \qquad\qquad+\tfrac{\tau^{2}}{8}\left[\e^{\frac{\tau}{2}b(u_{1})}b^{2}(u_{0})\phi_{3}\left(\tfrac{\tau}{2}b(u_{0})\right)+b^{2}\left(u_{1}\right)\phi_{3}\left(\tfrac{\tau}{2}b(u_{1})\right)\right]\\
 & \qquad=I+\tfrac{\tau}{8}\left[b(u_{0})+3b(u_{1})\right]\\
 & \qquad\qquad+\tfrac{\tau^{2}}{8}\left[b^{2}(u_{1})\phi_{2}\left(\tfrac{\tau}{2}b(u_{1})\right)+\tfrac{1}{2}b(u_{1})\phi_{1}\left(\tfrac{\tau}{2}b(u_{1})\right)b(u_{0})\right] \\ 
 &\qquad\qquad+\tfrac{\tau^{2}}{8}\left[ \e^{\frac{\tau}{2}b(u_{1})}b^{2}(u_{0})\phi_{3}\left(\tfrac{\tau}{2}b(u_{0})\right)+b^{2}\left(u_{1}\right)\phi_{3}\left(\tfrac{\tau}{2}b(u_{1})\right) \right].
\end{align*}
The term of first order vanishes and the term of second order is easily
shown to yield an additional order. In addition, we can easily bound
the remainder terms. 

Therefore, we have shown consistency for $S_{\tau}$. The extension of this result to $S_{\tau}^{(i)}$, for $i\geq2$, follows immediately from Theorem \ref{thm:fixedpoint}. 
\end{proof}
In the previous theorem we have established consistency and are thus
in a position to show convergence. Before stating this result let
us note that since we assume convergence of the (classic) Strang splitting
scheme $M_{\tau}$, it is not necessary to solve the significantly
more involved problem of nonlinear stability. Much effort has been
devoted to solve this problem for a variety of partial differential
equations found in the literature (see e.g.~\cite{time_analysis}
for the Vlasov--Poisson equations, \cite{lubich:2008} for the
Schrödinger--Poisson and cubic nonlinear Schrödinger equation, and \cite{holden2013} for the KdV equation). Therefore, we show in the following theorem that if convergence has
already been established for the (classic) Strang splitting scheme
this is also true for the iterated scheme we consider in this paper.
\begin{thm} \label{thm:convergence-hyperbolic}
(Convergence in the hyperbolic case). Suppose that assumption \eqref{eq:assumption_1} is satisfied and that $M_{\tau}$
is convergent of order two. Then it holds that (for $i\geq2$) 
\[
\left\Vert S_{\tau}^{(i)}u_{k}-u(k\tau)\right\Vert \leq C\tau^{2},
\]
where $C$ depends on solution $u(t)$ and the final time $T$.
\end{thm}
\begin{proof}
We once again proceed by comparing $S_{\tau}^{(i)}$ to $M_{\tau}$.
First, let us write 
\[
S_{\tau}^{(i)}u_{k}-M_{\tau}\overline{u}_{k}=S_{\tau}^{(i)}u_{k}-M_{\tau}u_{k}+M_{\tau}u_{k}-M_{\tau}\overline{u}_{k},
\]
where we have denoted the numerical approximation of the (classic)
Strang splitting scheme $M_{\tau}$ at time $k\tau$ by $\overline{u}_{k}$.
This immediately gives
\begin{align}
e_{k+1} & \leq\left\Vert S_{\tau}^{(i)}u_{k}-M_{\tau}u_{k}\right\Vert +\left\Vert M_{\tau}u_{k}-M_{\tau}\overline{u}_{k}\right\Vert \nonumber \\
 & \leq C\tau^{3}+\left\Vert M_{\tau}\right\Vert e_{k}.\label{eq:recursion-inequality}
\end{align}
In the equation above we have used $e_{k}$ to denote the difference 
to the (classic) Strang splitting scheme at step $k$, i.e.
\[
e_{k}=\left\Vert S_{\tau}^{(i)}u_{k}-M_{\tau}\overline{u}_{k}\right\Vert =\left\Vert \left(S_{\tau}^{(i)}\right)^{k}u_{0}-\left(M_{\tau}\right)^{k}u_{0}\right\Vert .
\]
Since by assumption $M_{\tau}$ is a stable numerical scheme we know
that 
\[
\sup_{0<k\tau\leq T}\Vert M_{\tau}^{k}\Vert\leq C.
\]
Now we know that $e_{0}=0$ and thus it can be easily shown that the
solution of the recursion inequality given in equation \eqref{eq:recursion-inequality}
satisfies $e_{k+1}\leq C\tau^{2}$,
which is the desired result.
\end{proof}

\section{Convergence for parabolic differential equations\label{sec:Convergence-for-parabolic}}

Let us note that Theorem \ref{thm:(Consistency)} and Theorem \ref{thm:convergence-hyperbolic} have an immediate extension to the parabolic case. However, due to the parabolic smoothing property we can weaken assumption \ref{eq:assumption_1} in the following sense: there exists a constant $C$ such that for all $0<t\leq T$
\begin{equation}
	\left\Vert b^{i}(\tilde{u})b^{(\ell)}(\tilde{u})b^{j}(\tilde{u})u(t)\right\Vert \leq \frac{C}{t^{\alpha}}
	,\qquad \ell+i+j\leq3,\label{eq:assumption_2}
\end{equation}
holds for $0\leq\alpha\leq1$.  Again $\tilde{u}$ represents
a generic function that is constructed in the manner described in section
\ref{sec:Convergence-for-hyperbolic}, but with the additional restriction
that 
$s\geq t$ holds
true (see the paragraph following assumption \eqref{eq:assumption_1}). 

The following theorem states that to achieve second order in time
(for $\alpha<1$) it is sufficient to assume less regularity in the
initial value than is required for the same convergence order in the
case hyperbolic differential equations.
\begin{thm}
(Convergence in the parabolic case). Suppose that assumptions \eqref{eq:assumption_2} is satisfied and that $M_{\tau}$ is
convergent of order two. Then it holds that
\[
\left\Vert S_{\tau}^{(i)}u_{k}-u(k\tau)\right\Vert \leq C\tau^{2}\log k.
\]

If assumption \eqref{eq:assumption_2} holds for $\alpha<1$, we get
a scheme of second order, i.e. 
\[
\left\Vert S_{\tau}^{(i)}u_{k}-u(k\tau)\right\Vert \leq C\tau^{2}.
\]
\end{thm}
\begin{proof}
We follow the proof presented in \cite{jahnke2000}. By using assumption
\eqref{eq:assumption_2} we can deduce from Theorem \ref{thm:(Consistency)}
that 
\[
e_{k+1}\leq C\frac{\tau^{3}}{k\tau}+\left\Vert M_{\tau}\right\Vert e_{k},
\]
where $C$ depends on the norm of the initial value. The above recursion
can be estimated to give
\[
e_{k+1}\leq C\tau^{2}s_{k}+\left\Vert M_{\tau}\right\Vert ^{k}e_{0},
\]
where 
\[
s_{k}=\sum_{i=1}^{k}\frac{1}{i}
\]
is the harmonic series. Since we know that $e_{0}=0$ we get the bound
\[
e_{k+1}\leq C\tau^{2}\log k,
\]
as desired.
\end{proof}

\section{High order composition methods\label{sec:high-order-composition}}

In the previous sections we have exclusively considered the iterated
Strang splitting scheme $S_{\tau}^{(i)}$. Now, we will
describe how this (almost) symmetric Strang splitting scheme can be
used to construct splitting methods of higher order. To that end we
will first discuss composition for a symmetric one-step method.

Suppose that a symmetric one-step method $\Phi_{\tau}$ is of even
order $p$. Then a method of order $p+2$ can be constructed by composition
(see e.g. \cite[p.~43]{hairer2006}). More specifically, under suitable
conditions on $\gamma_{1},\gamma_{2},\gamma_{3}$, the method $\Phi_{\gamma_{3}\tau}\circ\Phi_{\gamma_{2}\tau}\circ\Phi_{\gamma_{1}\tau}$
is of order $p+2$. Therefore, we are able to construct methods of
arbitrarily high even order. The cost, in terms of a single evaluation of
the corresponding second order method, is given by $3^{p/2-1}$. For
$p=4$, for example, the corresponding method is the well-known triple
jump scheme. 

That such an approach can be extended to the case of the almost symmetric
Strang splitting method $S_{\tau}^{(i)}$ has been shown, in the case
of non-stiff ordinary differential equations, in \cite{einkemmer2013}.
In this section we will extend the results presented in the before
mentioned paper to the time integration of partial differential equations.
Recall that if the following system of equations is satisfied
\begin{align}
 & \gamma_{1}+\gamma_{2}+\gamma_{3}=1,\label{eq:composition_gamma_1}\\
 & \gamma_{1}^{p+1}+\gamma_{2}^{p+1}+\gamma_{3}^{p+1}=0,\label{eq:composition_gamma_2}
\end{align}
the composition results in a scheme of order at least $p+1$. By a
symmetry argument we can then deduce that the order is indeed $p+2$.
This symmetry argument requires that $\Phi_{\tau}^{*}-\Phi_{\tau}=\mathcal{O}\left(\tau^{p+3}\right)$,
where we have used $\Phi_{\tau}^{*}$ to denote the adjoint method
of $\Phi_{\tau}$, and $\gamma_{1}=\gamma_{3}$ (see \cite{einkemmer2013}). 

The single real solution that simultaneously satisfies \eqref{eq:composition_gamma_1}
and \eqref{eq:composition_gamma_2} is given by
\[
\gamma_{1}=\gamma_{3}=\frac{1}{2-2^{1/(p+1)}},\quad\qquad\gamma_{2}=-2^{1/(p+1)}\gamma_{1}.
\]
Note that in order to perform this scheme we have to compute at least
one negative time step. In fact, this is shown to be necessary if
the desired order of the composition method is strictly larger than
two (see e.g. \textcolor{black}{\cite{blanes2005}}). While it is
possible to take negative time steps in the case of hyperbolic equations,
in parabolic problems this almost certainly leads to numerical instabilities
as, for example, roundoff errors are exponentially amplified in each
backward step. To remedy this situation in\textcolor{black}{{} \cite{hansen2009}}
and \cite{castella2009} it has been pointed out that the system of
equations given by \eqref{eq:composition_gamma_1} and \eqref{eq:composition_gamma_2}
admits complex solutions with positive real part. We will discuss
the implementational ramifications of using complex values in the
entire computation in section \ref{sec:brusselator}.

Now, let us turn to the main result of this section. That is, we will
indicate how the analysis conducted above can be extended from
the non-stiff ordinary differential equations considered in \cite{einkemmer2013}
to the partial differential equations considered in this paper. For such an analysis appropriate regularity assumptions are indispensable. In particular we have to control the remainder terms in the expansions.
Now, under the assumption that the
symmetric splitting scheme is consistent of order $p$, it holds that the local error of $\Phi_{\gamma_i \tau}$ satisfy
\begin{align*}
e_{1} & =C(u_0)\left(\gamma_{1}\tau\right)^{p+1}+\tau^{p+2}R_1(u_0)\\
e_{2} & =C(u(\gamma_1 \tau))\left(\gamma_{2}\tau\right)^{p+1}+\tau^{p+2}R_2(u(\gamma_1 \tau))\\
e_{3} & =C(u(\gamma_1+\gamma_2) \tau))\left(\gamma_{3}\tau\right)^{p+1}+\tau^{p+2}R_3(u(\gamma_1+\gamma_2)),
\end{align*}
where $R_{i}(\cdot)$ denotes the (bounded) remainder term. Note that
this, in fact, requires one more application of the $b$ operator
as computed in the previous section. However, it should be clear that
if the solution is sufficiently regular, such an expansion can be
accomplished by the straightforward, but tedious, extension of sections
\ref{sec:Convergence-for-hyperbolic} and \ref{sec:Convergence-for-parabolic}.
Then it is possible to generalize the proof given above by employing
Lemma \ref{lem:first-order-bound}. Let $E_{i}$ denote the full error, i.e.
\begin{align*}
	E_1 &= \Phi_{\gamma_1 \tau}(u_0) - u(\gamma_1 \tau) \\
	E_2 &= \Phi_{\gamma_2 \tau} \circ \Phi_{\gamma_1 \tau} (u_0) - u( (\gamma_1+\gamma_2) \tau) \\
	E_3 &= \Phi_{\gamma_3 \tau} \circ \Phi_{\gamma_2 \tau} \circ \Phi_{\gamma_1 \tau} (u_0) - u(\tau). 
\end{align*}
Then we have
\[
E_{i}=e_{i}+\tau\tilde{R}_{i},
\]
with a bounded remainder term $\tilde{R}_{i}$. The total error of
the composition method is then given by
\[
E_{1}+E_{2}+E_{3}=e_{1}+e_{2}+e_{3}+\tau^{p+2}Q,
\]
where $Q$ is a bounded remainder term. Thus, we get
\[
e_{1}+e_{2}+e_{3}=\tau^{p+1}\sum_{i=1}^{3}C\left(u_{i-1}\right)\gamma_{i}^{p+1}.
\]
To recover the result for the non-stiff case we have to expand $C(u_{1})$
and $C(u_{2})$. This can be done, for example, if the coefficient
$C(\cdot)$ is differentiable with respect to the (local) initial
value. All the estimates derived in sections \ref{sec:Convergence-for-hyperbolic}
and \ref{sec:Convergence-for-parabolic} have this property.

To conclude this section, let us note that the argument provided establishes
that, under suitable assumptions on the regularity of the solution,
the composition method is of order $p+1$. Due to the fact that $S_{\tau}^{(p+2)}$ is symmetric up to order $p+2$, the method is in fact of order $p+2$.


\section{Numerical simulations\label{sec:Numerical-simulations}}

Having established the theoretical convergence rates for the Strang
splitting scheme $S_{\tau}^{(i)}$ and its composition in the previous
sections, we now turn our attention to illustrate these results by conducting
numerical simulations. To that end we have chosen the well-known Brusselator
as a parabolic system. The Brusselator consists of two coupled advection-reaction
equations and provides a valuable prototype for many similar, but
often considerably more complicated, systems that model the interaction
of chemical species. As an example for a hyperbolic system we have
chosen the KdV equation in a single dimension. This equation is then
split in such a manner that the Burgers type nonlinearity can be solved
as a position-dependent but linear advection. For both of these examples
a \texttt{C++} program has been written that employs the \texttt{fftw}
library and is parallelized using \texttt{OpenMP}.

\subsection{A parabolic system: the Brusselator\label{sec:brusselator}}

As a example of a parabolic system we consider the Brusselator (see
e.g. \cite[p.~152]{hairer2}). This system is given by
\begin{align*}
\partial_{t}u & =\alpha\Delta u+\left(uv-\beta\right)u+\delta,\\
\partial_{t}v & =\alpha\Delta v-u^{2}v+\gamma u,
\end{align*}

where $u(t,x,y)$ and $v(t,x,y)$ are the two unknowns that usually
represent the concentration of a chemical species. We have chosen
the reaction parameters as $\beta=4.4$, $\gamma=3.4$, $\delta=1$,
and employ a relatively weak diffusion with $\alpha=10^{-2}$. The
equations above have to be supplemented by suitable initial conditions.
For the purpose of the numerical simulations conducted in this section
the following initial values are used
\begin{align*}
u_{0}(x,y) & =22y(1-y)^{3/2}\left(1+\cos(10\pi x)\right),\\
v_{0}(x,y) & =27x(1-x)^{3/2}\left(1+\sin(10\pi x)\right)
\end{align*}
and all computations are carried out on the domain $[0,1]^{2}$ with
periodic boundary conditions.

In the context of the splitting scheme outlined in section \ref{sec:Introduction},
let us define
\[
A\left[\begin{array}{c}
u\\
v
\end{array}\right]=\left[\begin{array}{c}
\alpha\Delta u\\
\alpha\Delta v
\end{array}\right].
\]
The corresponding flow can be computed efficiently by employing two
discrete Fourier transforms. In addition, the nonlinearity is
represented by (note that this choice is not unique)
\[
b\left(u_{\star},v_{\star}\right)=\left[\begin{array}{cc}
u_{\star}v_{\star}-\beta & 0\\
\gamma & -u_{\star}^{2}
\end{array}\right]
\]
and
\[
d=\left[\begin{array}{c}
\delta\\
0
\end{array}\right].
\]
Thus, we have to compute the solution of a linear system in two variables
which can be done analytically. Note that in this case we could, in
fact, find, by a rather tedious calculation, an analytical solution
for the complete nonlinearity. However, this is no longer possible
if one considers either additional variables, i.e. additional chemical
species, or alternatively a more complicated nonlinearity. However,
in such a case the splitting scheme outlined here would not suffer
any further difficulty or even a loss of efficiency as the nonlinearity
in question still reduces to a linear system that is easily solved
by standard methods.

\begin{figure}[h!]
\noindent \begin{centering}
\includegraphics[width=14cm]{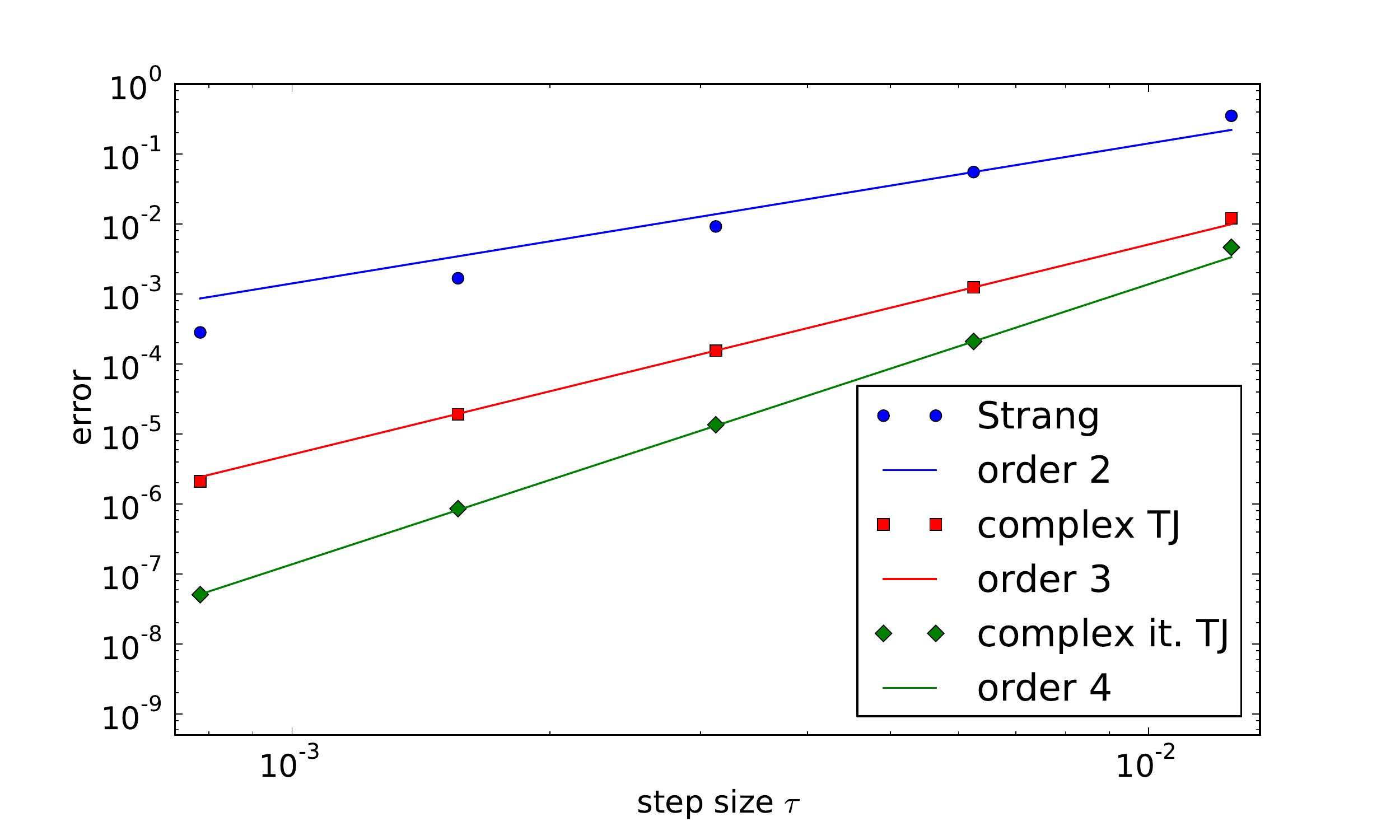}
\par\end{centering}

\caption{Order plot for the Brusselator integrated up to $T=0.25$ with 
	$2^{10}$ grid points per dimension (in total we have $2^{21}$ grid points). The error in time (in the discrete infinity
norm) is computed by comparing the numerical solution for a given
$\tau$ with a reference solution for which a sufficiently small time
step is chosen. For the triple jump and implicit triple jump scheme
complex arithmetics is employed. In addition, for each scheme a line
with slope equal to the expected order is displayed as well. \label{fig:brusselator_order}}
\end{figure}

In the numerical simulation we will compare three splitting schemes.
The (classic) Strang splitting scheme $M_{\tau}$ which is expected
to be of order two. The naive triple jump scheme which is constructed
by composition from $M_{\tau}$. Since $M_{\tau}$ is not symmetric
we expect this to be a third order scheme. Finally, we will consider
the triple jump scheme constructed by composition from $S_{\tau}^{(4)}$.
Since $S_{\tau}^{(4)}$ is symmetric up to order four we expect that
this approach results in a fourth order scheme. In fact, this is the
main motivation for our approach as it enables us to construct high
order (for example, fourth or sixth order methods) by employing the
well known composition rules. The numerical results shown in Figure
\ref{fig:brusselator_order} confirm the expected order for all numerical
schemes considered.

\begin{figure}[h!]
\noindent \begin{centering}
\includegraphics[width=14cm]{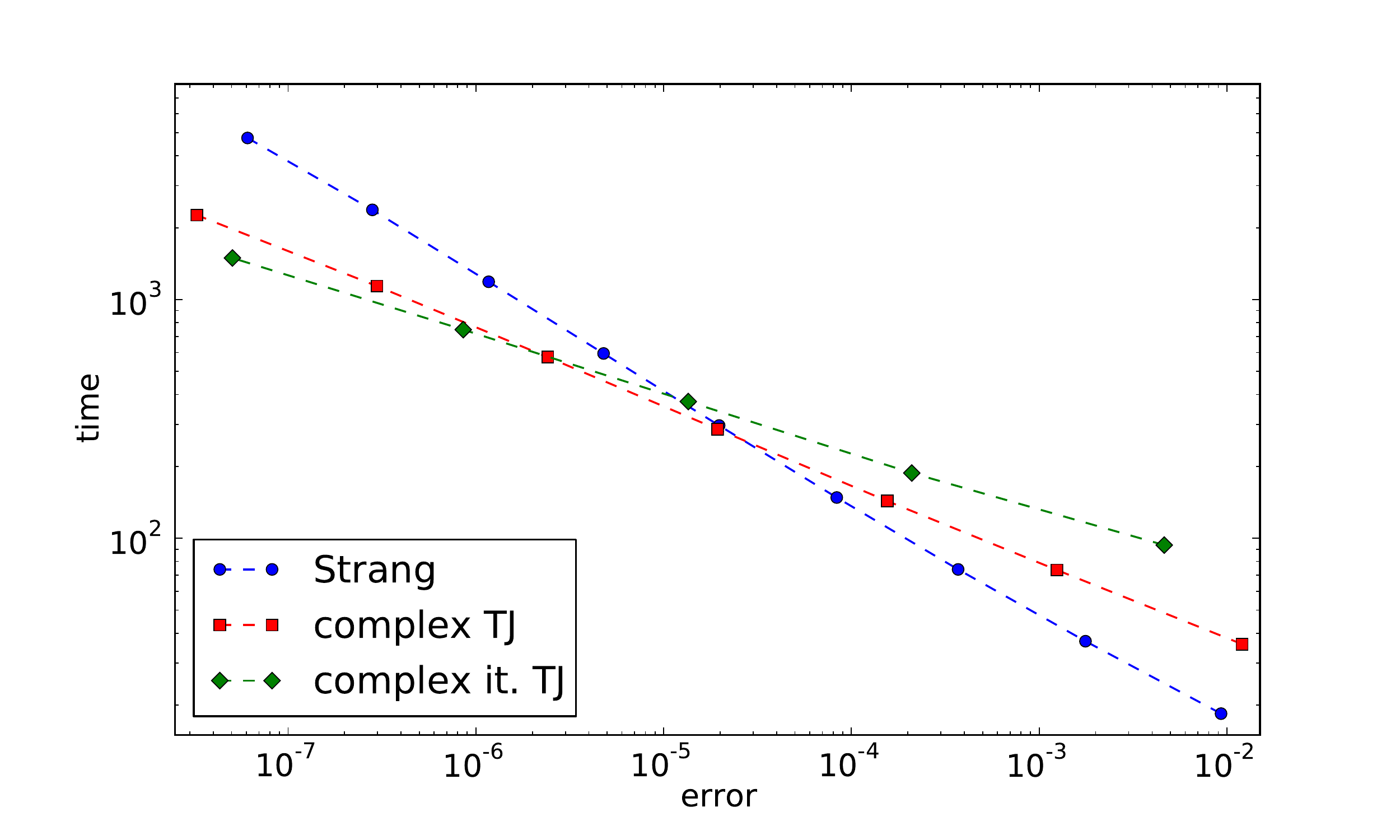}
\par\end{centering}

\caption{Work-precision plot for the Brusselator integrated up to $T=0.25$
	with $2^{10}$ grid points per dimension (in total we have $2^{21}$ grid points). The error in time (in the discrete
infinity norm) is computed by comparing the numerical solution for
a given $\tau$ with a reference solution for which a sufficiently
small time step is chosen. In addition, for each scheme a line with
slope equal to the inverse of the expected order is displayed as well.
\label{fig:Effort-plot-Brusselator}}
\end{figure}

From Figure \ref{fig:brusselator_order} we can clearly see that if
similar step sizes are taken the error made by the the fourth order
iterated triple jump scheme is significantly less than the error of
both the third order triple jump scheme as well as the second order
(classic) Strang splitting scheme. However, to discuss whether the
high order schemes constructed here can also provide a significant
gain in efficiency we have to plot the simulation time as a function
of the error. This is done in Figure \ref{fig:Effort-plot-Brusselator}.

It is shown that for high precision requirement (or equivalently long
integration times) the use of the fourth order iterated triple jump
scheme (with $4$ iterations, i.e. based on $S_{\tau}^{(4)}$) results
in a significant increase in efficiency. Also note that our fourth
order scheme is superior to the naive triple jump scheme for almost
any precision requirements. However, due to the overhead involved
in the use of complex arithmetics for small precision requirements
and short integration intervals the Strang splitting scheme is clearly
the preferred choice.

\subsection{A hyperbolic system: the KdV equation}

As an example of a hyperbolic system we consider the KdV (Korteweg--de
Vries) equation in a single dimension. It is given by
\[
\partial_{t}u(t,x)+\partial_{x}^{3}u(t,x)+u(t,x)\partial_{x}u(t,x)=0.
\]

It is shown in \cite{kametaka1969} and \cite{bona1975} that for
sufficiently regular initial values the regularity of the solution
to the KdV equation is not diminished as it is evolved in time (although
this does not rule out the appearance of high frequency oscillations).
For the purpose of this section we will consider two initial values
with extremely different dynamical behavior.

First, let us consider the following initial value
\begin{equation}
u_{0}(x)=u(0,x)=\frac{12}{\cosh^{2}x}.\label{eq:sechsoliton}
\end{equation}

The exact solution is a soliton that travels to the right with speed
$4$, i.e., the solution is given by (see e.g. \cite{klein2007numerical})
\[
u(t,x)=u_{0}(x-4t).
\]

For our numerical studies we consider the domain $[-20,20]$ in order
to limit artifacts which originate from the fact that periodic boundary
conditions are imposed on the domain of finite length that is used in the
numerical simulations. We integrate up to a final time $T=0.4$.

Second, we consider an initial value where oscillations appear. The
so called Schwartzian initial value is given by (see e.g. \cite{klein2011})
\begin{equation}
u_{0}(x)=\frac{12x\tanh\vert x\vert}{\vert x\vert\cosh^{2}x}.\label{eq:schwartzian}
\end{equation}

In this instance we integrate only up to $T=0.05$ and employ the
domain $[-4\pi,4\pi]$ with periodic boundary conditions (following
the same argument as given above).

Now let us turn our attention to the splitting scheme employed. In
the framework of section \ref{sec:Introduction} the linear problem
is defined as 
\[
A=-\partial_{x}^{3}
\]

and can be solved efficiently by employing two discrete Fourier transforms.
Instead of the full nonlinearity we numerically compute the evolution
corresponding to
\begin{equation}
b\left(u_{\star}\right)=u_{\star}(x)\partial_{x},\qquad d=0.\label{eq:KdV-linearized}
\end{equation}

This is significantly less involved than solving the full Burgers
type nonlinearity. In fact, \cite{klein2011} states that solving
the full nonlinearity invalidates splitting as a viable approach to
numerically solve the KP (Kadomtsev--Petviashvili) equation (note
that, compared to the KdV equation, the two dimensional KP equation
includes an additional diffusive term in the additional variable).
However, the linear problem given by equation \eqref{eq:KdV-linearized}
can be computed numerically to high precision, for example, by using
the (exact) exponential of a finite difference stencil (as it is often
done in the context of exponential integrators, for example). We
intended to employ the \texttt{Expokit }package to compute the exponential
of the a seven-point stencil that is used to approximate the first
derivative. However, in case of the soliton solution the results were
not satisfactory except for very small step sizes. Therefore, we used
the \texttt{SUNDIALS CVODE} solver, with a prescribed tolerance of
$10^{-12}$, for this example%
\footnote{Unfortunately, to our knowledge, there are no packages that are written
in \texttt{C++} available to compute the matrix exponential. We have
tried both \texttt{Expokit} as outlined above and the (unsupported)
matrix exponential provided by the \texttt{SPARSEKIT} package (using
our own \texttt{Fortran} to \texttt{C} bindings). Even though the
source code of \texttt{SPARSEKIT} is much more readable, it is only
able to compute an approximation up to a tolerance of $10^{-6}$.
This, is even true for the diagonal example provided as part of the
package. None of these packages are parallelized. However, let us
mention here that there are a number of \texttt{Python} and \texttt{Matlab}
implementations.%
}.

\begin{figure}[h!]
\noindent \begin{centering}
\includegraphics[width=14cm]{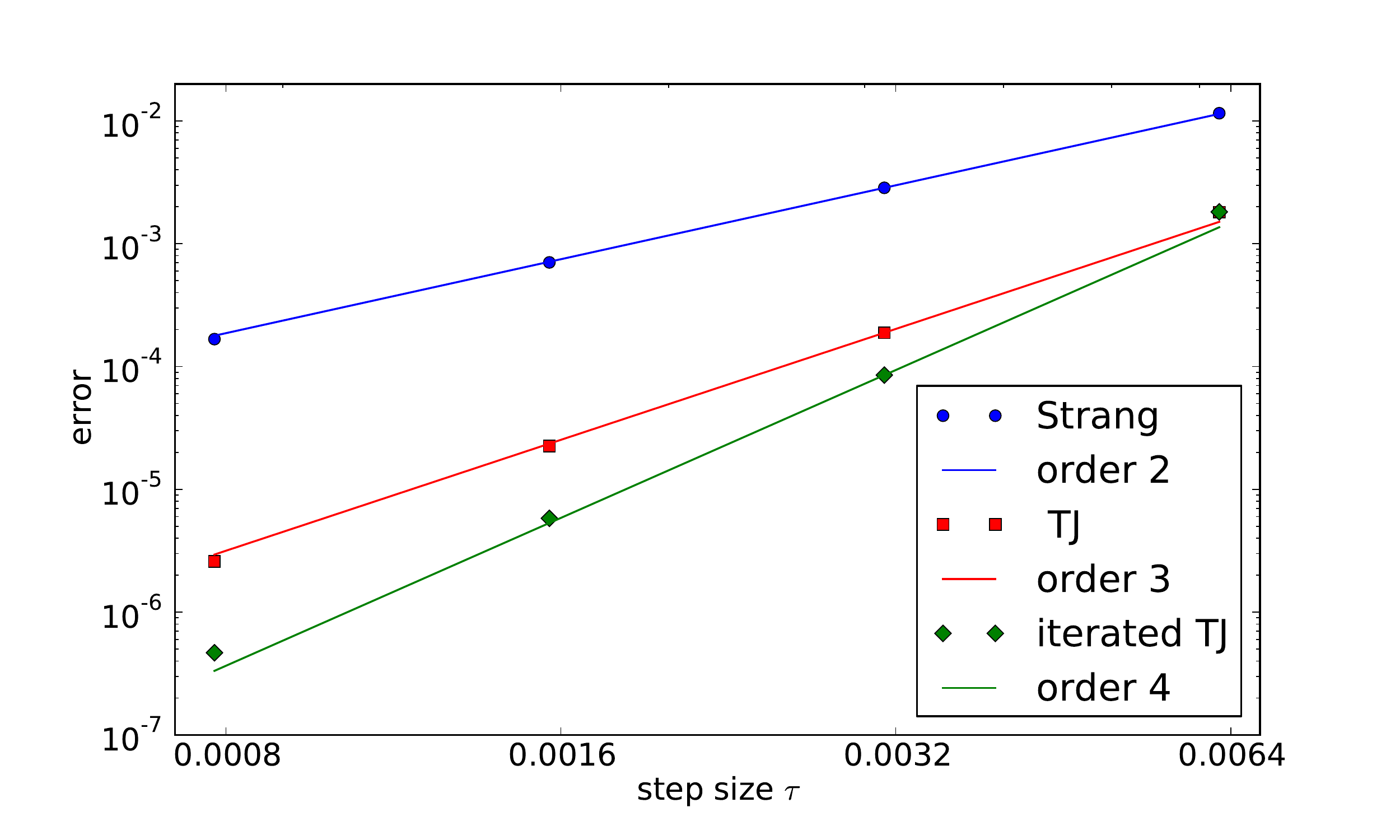}
\par\end{centering}

\caption{Order plot for the KdV equation integrated up to $T=0.4$ with $2^{10}$
grid points for the soliton initial value given in \eqref{eq:sechsoliton}.
The error in time (in the discrete $L^{2}$ norm) is computed by comparing
the numerical solution for a given $\tau$ with a reference solution
for which a sufficiently small time step has been chosen. In addition,
for each scheme a line with slope equal to the expected order is displayed
as well. \label{fig:order-plot-KdV-sechsoliton}}
\end{figure}

\begin{figure}[h!]
\noindent \centering{}\includegraphics[width=14cm]{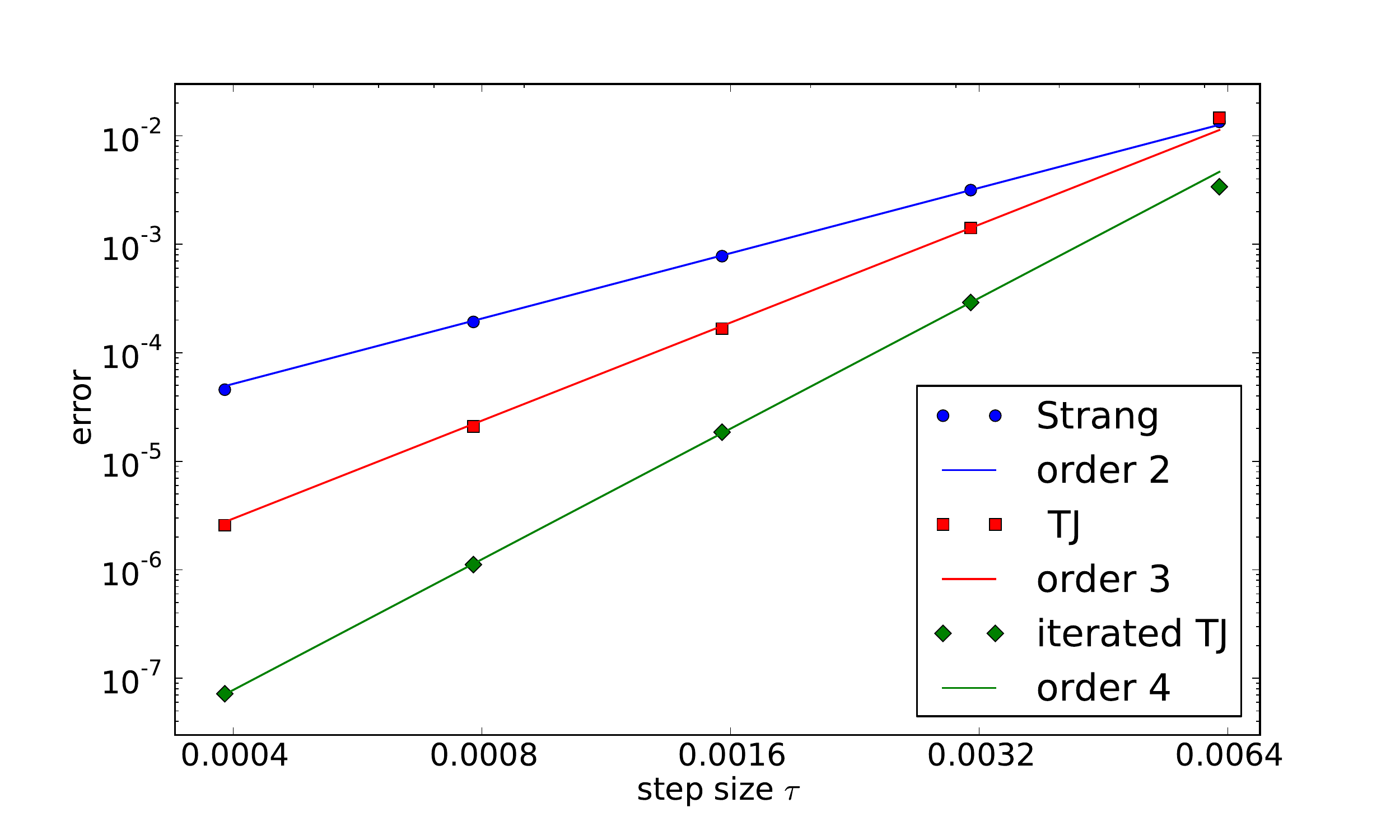}\caption{Order plot for the KdV equation integrated up to $T=0.05$ with $2^{11}$
grid points for the Schwartzian initial value given in \eqref{eq:schwartzian}.
The error in time (in the discrete $L^{2}$ norm) is computed by comparing
the numerical solution for a given $\tau$ with a reference solution
for which a sufficiently small time step has been chosen. In addition,
for each scheme a line with slope equal to the expected order is displayed
as well. \label{fig:order-plot-KdV}}
\end{figure}

The numerical results are shown in Figure \ref{fig:order-plot-KdV-sechsoliton}
(for the soliton initial value given in equation \eqref{eq:sechsoliton})
and in Figure \ref{fig:order-plot-KdV} (for the Schwartzian initial
value given in equation \eqref{eq:schwartzian}). Note that the numerical simulation matches the predicted order of
the schemes studied very well. Thus, we conclude that the data obtained
are clearly consistent with the analysis conducted in sections \ref{sec:Convergence-for-hyperbolic}
and \ref{sec:high-order-composition}.

\section{Conclusion and Outlook}

We have rigorously shown that the proposed iterated Strang splitting
scheme is of second order in time and due to its symmetry properties
can be used to construct methods of arbitrary (even) order by composition.
The main assumption we have made is that the classic Strang splitting
method is convergent of order two (i.e., that we deal with a problem
for which applying a splitting scheme of second order is sensible).
Further, a technical assumption on the nonlinearity is made. This
assumption reduces, in the linear limit, to the statement that we
have to bound an appropriate number of application of the operator
in question to the exact solution (assumptions of that type have been
used in much of the literature to show convergence of splitting methods
for linear partial differential equations). 

In addition, we have provided an argument demonstrating that our iterated
Strang splitting can be used, similar to the case of ordinary differential
equations, to construct higher order methods by composition if sufficient
regularity of the exact solution can be assumed. This has been verified
up to order four for both the Brusselator system (a parabolic problem)
and the KdV equation (a hyperbolic partial differential equation).
In both instances the iterated fourth order method is shown to provide
superior performance, in case of medium to high precision requirements
(or equivalently long integration times), compared to both the (classic)
Strang splitting scheme as well as the (classic) triple jump scheme
(which is a method of order three). We also conclude that the necessity
of using complex precision arithmetics, in the case of parabolic problems,
does not negate the performance gain we expect from higher order methods.

In this paper we have only considered methods of order up to four.
However, composition can be used to construct methods of arbitrary
(even) order. To provide a clear picture of the efficiency gain expected
for such high order methods, especially in the context of the semi-Lagrangian
approach discussed in this paper, longer integration times as well
as a space discretization with significantly more grid points has
to be considered. Since we have not considered parallelization and
other computing aspects in this paper, we will consider such an implementation
as a subject of further research. 

\bibliographystyle{plain}
\bibliography{pde-iterStrang}

\end{document}